\documentclass[12pt]{article}
\usepackage{amssymb}
\newtheorem{proposition}{Proposition}[section]  %proposition
\newcommand{\bprop}{\medskip\begin{proposition}  \it}
\newcommand{\eprop}{\end{proposition} \hfill  \\}

\newcommand{\bed}{\begin{displaymath}}
\newcommand{\eed}{\end{displaymath}}

\newcommand{\ba}{\begin{array}}
\newcommand{\ea}{\end{array}}
\newcommand{\beq}{\begin{equation}}
\newcommand{\eeq}{\end{equation}}
\newcommand{\be}{\begin{equation}}
\newcommand{\ee}{\end{equation}}
\newcommand{\bea}{\begin{eqnarray}}
\newcommand{\eea}{\end{eqnarray}}

\def\id{\mbox{$\mathop{\mbox{\rm id}}$}}

\def\ot{{\otimes}}
\def\Z{{\mathbb Z}}
\def\N{{\mathbb N}}
\def\C{{\mathbb C}}

\usepackage{amsmath}
\usepackage{amsfonts}

\newcommand{\mZ}{\mathbb{Z}}

\newcommand{\mdeg}{\text{deg}}
\newcommand{\mId}{\text{id}}
\def\sw#1{{\sb{(#1)}}}
\def\eps{\varepsilon}
\def\can{\text{can}}

\newcommand{\mAd}{\text{Ad}}
\newcommand{\mD}{\text{d}}
\begin{document}
\baselineskip17.5pt
\setcounter{page}{0}
\thispagestyle{empty}
%%%%%%%%%%%%%%%%%%%%%%%%%%%%%%%%%%%%%%%%%
%\begin{flushright}
%{\large\bf Preliminary version}
%\end{flushright}
%%%%%%%%%%%%%%%%%%%%%%%%%%%%%%%%%%%%%%%%%
%
%~\vspace{1.5cm}
%
\begin{center}{
\LARGE \bf
Hopf fibration and monopole connection \\
over the contact quantum spheres.
}
\end{center}
\vspace{1cm}
\centerline{\Large
Tomasz\ Brzezi\'nski}
\vspace{-5mm}
\begin{center}
{\it Department of Mathematics, University of Wales Swansea\\
   Singleton Park, Swansea SA2 8PP, U.K.;}
~T.Brzezinski@swansea.ac.uk
\end{center}
\vspace{5mm}
\centerline{\Large
Ludwik\ D\c abrowski}
\vspace{-5mm}
\begin{center}
{\it Scuola Internazionale Superiore di Studi Avanzati,\\
Via Beirut 2-4, I-34014, Trieste, Italy;}
~dabrow@sissa.it
\end{center}
\vspace{5mm}
\centerline{\Large
Bartosz\ Zieli\'nski}
\vspace{-5mm}
\begin{center}
{\it Department of Mathematics, University of Wales Swansea,\\ Singleton
Park, Swansea SA2 8PP, U.K.;}
~mabpz@swan.ac.uk
\\ and \\
{\it Department of Theoretical Physics II, University of \L{}\'od\'z \\
Pomorska 149/153, 90-236 \L{}\'od\'z, Poland.}
\end{center}
\vspace{1cm}\noindent
\begin{abstract}\noindent
Noncommutative geometry of  quantised contact spheres introduced by
Omori-Maeda-Miyazaki-Yoshioka in \cite{OMMY,OMYM}
is studied. In particular it is proven that these spheres form a
noncommutative Hopf fibration in the sense of Hopf-Galois extensions.
The monopole (strong) connection is constructed,
and projectors describing projective
modules of all monopole charges are computed.
   \end{abstract}

\vspace{1.5cm}
\centerline{keywords:
Noncommutative geometry, quantum spheres, Hopf fibration, monopole}
\centerline{mathclass: Primary 58B34; Secondary 17B37.}

\section{Introduction}
Since its inception in \cite{BrzMaj:gau}, the geometric approach to
quantum group gauge theories based on quantum principal bundles or
Hopf-Galois extensions, attracted a considerable interest. The general
structure of quantum principal bundles and their generalisations known
as coalgebra principal bundles or coalgebra-Galois extensions
\cite{BrzMaj:coa,BrzMaj:geo} is well-known by now. In particular those
Hopf- and coalgebra-Galois extensions which admit a special kind of
connections, known as strong connections introduced in \cite{Haj:str}
seem to be  of special interest for non-commutative geometry a la
Connes. As recently revealed in  \cite{HajMaj:pro}, \cite{DabGro:str},
\cite{BrzHaj:rel} the modules of sections of quantum vector bundles
associated to such extensions are projective modules, i.e., vector
bundles of non-commutative geometry. This makes Galois-type extensions
with strong connections perfect candidates for (algebraic)
non-commutative principal bundles. For an extensive review of algebraic
background and recent progress in this field we refer to
\cite{BrzHaj:gal}.

Despite the fact that the structure of Galois-type extensions is very
well-known, it is  quite difficult to construct concrete examples of
such extensions which in addition admit strong connections. These are
technical rather than structural difficulties, hence they do not
undermine usefulness of the general theory. On the other hand interest
in any abstract theory is fuelled by non-trivial examples. For many
years essentially one non-trivial example of a quantum principal bundle
with strong connection was known. This was the original example of a
quantum Hopf fibration in \cite{BrzMaj:gau}. Recently, however, a number
of new non-trivial examples have been constructed. These include
examples obtained by patching of trivial quantum principal bundles
\cite{CalMat:con} as well as examples motivated by recent interest in
and several constructions of quantum and non-commutative spheres
initiated in \cite{ConLan:man} and  \cite{DabLan:ins} (see
\cite{Dab:gar} for a concise review of low dimensional cases).
An earlier
deformation of the 3-sphere introduced by Matsumoto in
\cite{Mat:sph} has been shown to give rise to a Hopf-Galois extension of
the (commutative) algebra of functions on the 2-sphere
\cite{BrzSit:mat}. Other examples, in particular \cite{ConDub:sph} and
\cite{BonCic:ins} seem also to fit perfectly into the scheme of
Galois-type extensions.

The aim of the present note is to reveal that quantum spheres
constructed in \cite{OMMY,OMYM}  as quantisations of contact structures
on spherical manifolds, and thus termed {\em contact quantum spheres}
give rise to quantum principal bundles with strong connections. These
are all non-standard deformations of the Hopf fibration. The resulting
connections appear to be deformations of the Dirac magnetic monopole
field. The general theory guides us then to  projectors for the
associated quantum line bundles with an arbitrary monopole charge.

   The paper is ogranised as follows. In Section~2 we describe the
non-standard deformation of the 3-sphere, obtained in \cite{OMMY} by the
quantisation of the contact structure on the 3-sphere. We show that
this sphere admits an action of the group $U(1)$. In algebraic terms it
means that there is a coaction of the algebra of Laurent polynomials in
one variable. Next we show that the fixed points of these action
(coaction)
describe the deformed algebra of the noncommutative 2-sphere, introduced
in
\cite{OMMY}. Finally we show that this structure defines
a quantum Hopf fibration (cf.\ \cite[Section~4]{OMYM}), namely that
the deformed algebra of the noncommutative 3-sphere is a Hopf-Galois
extension of the deformed algebra of the noncommutative 2-sphere.
In Section~3 we construct explicitly a strong connection on this
Hopf-Galois extension. This immediately gives rise to projectors for
quantum
associated line bundles, which we explicitly compute.
In Section~4 we address some questions about the problem of specifying
the formal
paramter $\mu$ to a numerical value and about the possibility of a
$C^*$-algebraic
version.

We work over a field of complex numbers $\mathbb{C}$ and unadorned
tensor product is over $\mathbb{C}$.

\section{Quantum Hopf fibration}

The algebraic structure underlying quantum group principal bundles is
provided by Hopf-Galois extensions. To construct such an extension one
needs the following data. First, let $H$ be a  Hopf algebra with
coproduct $\Delta:H\to H\otimes H$, counit $\eps: H\to \C$ and antipode
$S:H\to H$. Second, let $P$ be a right $H$-comodule algebra, i.e., an
algebra and a right $H$-comodule such that the coaction $\Delta_R:P\to
P\otimes H$ is an algebra map. For the coaction we use the Sweedler
notation
$\Delta_R(p) = p\sw 0\otimes p\sw 1$ (summation understood). With these
data one defines the {\em coinvariant} subalgebra of $P$ by
$$
B=\{x\in P\; |\; \Delta_{R}(x)=x\ot 1\}.
$$
Since $B$ is a subalgebra of $P$ there is an obvious inclusion map,
hence an extension, $B\to P$. Furthermore, $P$ is a $(B,B)$-bimodule,
hence one can consider the tensor product $P\otimes_B P$. The extension
$B\to P$ is called a {\em Hopf-Galois extension}, provided the canonical
left $P$-module, right $H$-comodule map
$$
\can :=(m_{P}\otimes\mId)\circ(\mId\otimes_{B}\Delta_{R}) : P\otimes_B
P\to P\otimes H,
$$
where $m_P$ denotes the product in $P$, is bijective. Explicitly, the
map $\can$ reads
$$
\can (x\otimes_{B} y) = xy\sw 0\otimes y\sw 1.
$$
The aim of this section is to show that the quantum contact 3-sphere
defined in
\cite{OMMY} is a Hopf-Galois extension of its algebra of coinvariants.
The latter is the algebra of functions of the quantum 2-sphere also
defined in \cite{OMMY}.

The polynomial $*$-algebra $A(S^3_\mu)$ that underlies the
Omori-Maeda-Miyazaki-Yoshioka quantum contact 3-sphere
is generated by a selfadjoint
$\mu$ and by $a, b, a^*, b^*$ with relations
\begin{subequations}
      \label{eq1}
\begin{gather}
	ba=ab,\qquad	ab^{\ast}=(1-\mu)b^{\ast}a,\label{eq1a}\\
	\mu a-a\mu=\mu a\mu, \quad \mu b-b\mu=\mu b\mu, \label{eq1c}\\
	aa^{\ast}-(1-\mu)a^{\ast}a=\mu, \quad bb^{\ast}-(1-\mu)b^{\ast}b=\mu,\\
	a^{\ast}a+b^{\ast}b=1.\label{eq1e}
\end{gather}
\end{subequations}
  From the above relations it follows that
\begin{equation}
      \label{eq1f}
      aa^* + bb^*=1+\mu.
\end{equation}
Note also that (\ref{eq1c}) is equivalent to
$$
\mu a(1+k\mu) = (1+(k+1)\mu) a\mu, \qquad \mu b(1+k\mu)=(1+(k+1)\mu)
b\mu ,
$$
for all $k\in \mZ$. We denote by $A'(S^3_\mu)$ the
polynomial $*$-algebra $A(S^3_\mu)$ with the generator $\mu$
required to be invertible (i.e. with adjoined $\mu^{-1})$.

It turns out convenient for our purposes to employ
a certain {\em $\mu$-regulated smooth algebra} $A^\infty(S^3_\mu)$,
which is defined and studied in \cite{OMMY}.
The algebra $A^\infty(S^3_\mu)$,
called {\em noncommutative contact algebra} on $S^3$,
contains densely the polynomial $*$-algebra $A'(S^3_\mu)$.
Also, it contains $f(\mu )$ for any formal power series $f$,
and the following relations are fulfilled
\begin{subequations}
      \label{series}
\begin{gather}
      a f(\mu )=f(\frac{\mu}{1+\mu} )a,\qquad   f(\mu ) a =
af(\frac{\mu}{1-\mu}),\\
      b f(\mu )=f(\frac{\mu}{1+\mu} )b, \qquad    f(\mu ) b =
bf(\frac{\mu}{1-\mu}).
\end{gather}
\end{subequations}
In particular, for all $k\in\mZ$ the elements $1+k\mu$
are invertible and have square root in $A^\infty(S^3_\mu)$.
In the sequel we shall need their inverses as well as the square roots,
which satisfy the following relations $\forall k\in \mZ$,
\begin{equation}
      \label{invs}
      a{\mu}({1+k\mu})^{-1}={\mu}({1+(k+1)\mu})^{-1}a,\qquad
b{\mu}({1+k\mu})^{-1}={\mu}({1+(k+1)\mu})^{-1}b
\end{equation}
and
\begin{equation}
      \label{sqrts}
a\sqrt{1+k\mu} = \frac{\sqrt{1 + (k+1)\mu}}{\sqrt{1+\mu}} a, \qquad
b\sqrt{1+k\mu} = \frac{\sqrt{1 + (k+1)\mu}}{\sqrt{1+\mu}} b.
\end{equation}

Although $\mu$ is a generator
it can be regarded  as a noncentral formal parameter,
cf.\ \cite{OMMY} for a precise meaning of this statement.
Also, from the defining relations of (\ref{eq1}) it is apparent that
$A^\infty(S^3_\mu)$
is a $\mZ$-graded algebra with the grading defined by
setting
$$\mdeg(a)=\mdeg(b)=1, \mdeg(a^{\ast})=\mdeg(b^{\ast})=-1, \mdeg(\mu)=0.$$
This in turn allows us to view $A^\infty(S^3_\mu)$
as a comodule algebra of the Hopf algebra $H$ of functions on $U(1)$.
Explicitly,
$H=\C[u,u^{-1}]$ is an algebra of Laurent polynomials in one variable
$u$ (i.e., $u^{-1}$ is the multiplicative inverse of $u$), and it is a
$*$-algebra with $u^* = u ^{-1}$.
A  Hopf algebra structure of $H$ is determined by
$\Delta (u) = u\ot u$, $\eps (u) = 1$ and $S(u) =  u^{-1}$.
The grading of $A^\infty(S^3_\mu)$ makes it a right comodule algebra
with the
coaction given on homogeneous elements by
$$
\Delta_R(x) = x\otimes u^{\mdeg (x)}.
$$
Thus explicitly on generators the coaction comes out as
$$\Delta_R (a) = a\ot u, \Delta_R (b) = b\ot u,
\Delta_R (a^*) = a^*\ot u^{-1}, \Delta_P (b^*) = b^*\ot u^{-1},
\Delta_P (\mu ) = \mu\ot 1.$$
Note that the coaction $\Delta_R$ is
compatible with the $*$-structure (it is a $*$-algebra homomorphism).

The definition of the coaction in terms of the grading immediately
implies that the coinvariant subalgebra coincides with the zero-degree
subalgebra, i.e.,
\begin{equation*}
A^\infty(S^2_\mu)
:=\{x\in A^\infty(S^3_\mu)\; |\;
\Delta_{R}(x)=x\ot 1\}=\{x\in A^\infty(S^3_\mu)\; |\; \mdeg(x)=0\}.
\end{equation*}
Using equation (\ref{invs}) it is immediate to verify
that $\mu$ is central element in $A^\infty(S^2_\mu) $.
It can be also seen that $A^\infty(S^2_\mu) $ is the commutant
of $\mu$ in $A^\infty(S^3_\mu) $.

The relations  (\ref{eq1}), (\ref{invs}) and (\ref{sqrts})
provide us with a deeper insight into the structure of
$A^\infty(S^2_\mu)$.
With their help we can establish that $A^\infty(S^2_\mu) $ contains a
(dense)
polynomial $*$-algebra $A'(S^2_\mu)$ generated by
$$X = X^* = a a^* -(\mu +1)/2 , \quad Z = a b^* , \quad  Z^* = b a^*,
$$
and self-adjoint (invertible) element $\mu$.
Notice that $A'(S^2_\mu)$ is contained strictly in the commutant
of $\mu$ in $A'(S^3_\mu) $, which coincides
with the grade zero subalgebra of $A'(S^3_\mu) $ and also
with the $H$-coinvariant subalgebra of $A'(S^3_\mu) $.
Similar observations hold for $A(S^2_\mu)$, defined as the $*$-algebra
obtained from $A'(S^2_\mu)$ by omitting the invertibility of $\mu$.

The relations in $A^\infty(S^2_\mu) $ are derived from the relations in
$A^\infty(S^3_\mu)$ and come out as
\begin{subequations}
      \label{eq3}
      \begin{gather}
	\mu X - X\mu =0,\; 	\mu Z - Z\mu =0,\label{eq3aa}\\
	XZ-ZX=-\mu Z,\label{eq3a}\\
	ZZ^{\ast}-Z^{\ast}Z=-2\mu X,\label{eq3b}\\
	\left(X+\frac{\mu}{2}\right)^{2}+ZZ^{\ast}=\frac{1}{4},\label{eq3c}
      \end{gather}
\end{subequations}
Note that from the relations above it follows that there is also a
second radial relation
$$
(X-\frac{\mu}{2})^2 + Z^*Z = \frac{1}{4}.
$$
Since $\mu$ is central, it could be possible to consider $\mu$
as a formal parameter and specify it to a numerical value,
we shall comment on this issue in the final section.

We claim that in a `dual' sense the quantum contact three sphere
$S^3_\mu$
is a total space of  a quantum $U(1)$-principal bundle
over the quantum two-sphere $S^2_\mu$, i.e.,
$A^\infty(S^3_\mu)$ is a Hopf-Galois extension of $A^\infty(S^2_\mu)$
with the structure Hopf algebra $H$.
This claim is proven by explicit construction
of the inverse to the canonical map.

To relieve the notation we write $P$ for $A^\infty(S^3_\mu)$
and $B$ for $A^\infty(S^2_\mu)$.
Consider the map
$\can^{-1}:P\otimes H\rightarrow P\otimes_{B}P$ defined
for all $x\in P$ and $n\in \N$ ~by
\begin{subequations}
      \label{eq4}
      \begin{align}
	\can^{-1}(x\ot u^{n})&=\sum_{k=0}^{n}{n\choose k}
         x(a^{\ast})^{n-k}(b^{\ast})^{k}
	\otimes_{B}b^{k}a^{n-k},\label{eq4a}\\
	\can^{-1}(x\ot u^{-n})&=x({1+n\mu})^{-1}\sum_{k=0}^{n}{n\choose k}
	a^{n-k}b^{k}\otimes_{B}(b^{\ast})^{k}(a^{\ast})^{n-k},\label{eq4b}
      \end{align}
\end{subequations}
where ${n\choose k}$ are the usual binomial coefficients.
Directly from the definition it follows
that  $\can^{-1}$  is a left $P$-module map. Furthermore, the degree
counting on the right hand side and the comparison of the powers of
$u$ immediately confirm that $\can^{-1}$ is a right $H$-comodule map.
Before we  prove that $\can^{-1}$ is the inverse map to $\can$ we note
that for all $n\in \N$,
\begin{subequations}
      \label{eq7}
      \begin{align}
	\sum_{k=0}^{n}{n\choose k}(a^{\ast})^{n-k}(b^{\ast})^{k}
	b^{k}a^{n-k}&=1,\label{eq7a}\\
	\sum_{k=0}^{n}{n\choose k}a^{n-k}b^{k}(b^{\ast})^{k}(a^{\ast})^{n-k}
	&=1+n\mu.
	\label{eq7b}
      \end{align}
\end{subequations}
The formulae (\ref{eq7}) are most easily proven by induction.
They are clearly satisfied for $n=1$. Next, assume that they hold for $n-1$
with $n\geq 2$.
Using equations (\ref{eq1e}), (\ref{eq1f}) and the well-known formula
\begin{equation*}
      \sum_{l=0}^{k}(-1)^{l}{n\choose{k-l}}={{n-1}\choose k},
\end{equation*}
observe that
\begin{subequations}
      \label{eq8}
      \begin{gather}
	\sum_{k=0}^{n}{n\choose k}(a^{\ast})^{n-k}(b^{\ast})^{k}
	b^{k}a^{n-k}=
	\sum_{k=0}^{n-1}{n-1\choose k}(a^{\ast})^{n-1-k}(b^{\ast})^{k}
	b^{k}a^{n-1-k},
	\label{eq8a}\\
	\sum_{k=0}^{n}{n\choose k}a^{n-k}b^{k}(b^{\ast})^{k}(a^{\ast})^{n-k}
	=
	\sum_{k=0}^{n-1}{n-1\choose k}a^{n-1-k}b^{k}(1+\mu)
	(b^{\ast})^{k}(a^{\ast})^{n-1-k}.
	\label{eq8b}
      \end{gather}
\end{subequations}
Then using (\ref{invs}) we conclude that equations (\ref{eq7a},
\ref{eq7b}) hold for all $n$.
Now we are in position to prove that $\can^{-1}$ is the inverse of
$\can$. Take $x\in P$ and $n\in \N$. Then
\begin{eqnarray*}
\can(\can^{-1}(x\ot u^n)) & = & \can \left(\sum_{k=0}^{n}{n\choose
k}x(a^{\ast})^{n-k}(b^{\ast})^{k}
	\otimes_{B}b^{k}a^{n-k}\right) \\
&=& \sum_{k=0}^{n}{n\choose k}x(a^{\ast})^{n-k}(b^{\ast})^{k}b^{k}a^{n-
k}\ot u^n = x\ot u^n,
\end{eqnarray*}
where the last equality follows from (\ref{eq7a}). Similarly, the use of
(\ref{eq7b}) confirms that
$$\can(\can^{-1}(x\ot u^{-n})) =x\ot u^{-n}.
$$
Conversely we need to check the equality $\can^{-1}(\can(x\ot_B y)) =
x\ot_B y$ for all $x,y\in P$. Since $P$ is a $\mZ$-graded algebra
suffices it to take homogeneous  $y$ of degree $n$. Suppose $n\geq 0$.
Then
$$
\can^{-1}((\can(x\ot_B y)) = \can^{-1}(xy\ot u^n) =
\sum_{k=0}^{n}{n\choose k}xy(a^{\ast})^{n-k}(b^{\ast})^{k}
	\otimes_{B}b^{k}a^{n-k}.
$$
Since $\mdeg(y) =n$, each of the $y(a^{\ast})^{n-k}(b^{\ast})^{k}$ has
degree 0, hence it is in  $B$ and we can write
$$
\can^{-1}((\can(x\ot_B y)) = \sum_{k=0}^{n}{n\choose k}x
	\otimes_{B}y(a^{\ast})^{n-k}(b^{\ast})^{k}b^{k}a^{n-k} = x\ot_B y,
$$
by (\ref{eq7a}). In the case of homogeneous $y$ of negative degree,
we use equation (\ref{eq7b}) to obtain the assertion. Thus we have proven
that $\can$ is a bijective map, i.e., $A^\infty(S^3_\mu)$ is a
Hopf-Galois
extension of $A^\infty(S^2_\mu)$ as claimed.

\section{Monopole connection and projectors of charge $n$}

\subsection{Strong connection}
Connections in quantum principal bundles are defined as colinear
splittings of one forms into horizontal and vertical parts (cf.\
\cite{BrzMaj:gau}). From the non-commutative geometry point of view a
special class of connections, introduced in  \cite{Haj:str} and known as
strong connections, is of particular interest.

Algebraically, one-forms on an algebra $P$ are defined as a $P$-bimodule
$\Omega^1(P)$ together with a linear map $\mD:P\to \Omega^1(P)$ which
satisfies the Leibniz rule $\mD (pp') = \mD (p)p' + p\mD (p')$. In what
follows we restrict ourselves to the particular example of differential
forms, known as universal differential forms and denoted by
$\Omega^1P$.   The bimodule $\Omega^1P$ is defined as the kernel of the
multiplication map $m_P$, and the differential in this case is given by
the map
$$
\mD  : P\to \Omega^1P = \ker m_P, \qquad p\mapsto 1\otimes p - p\otimes
1.
$$
   In the universal differential calculus case, strong connections  on a
Hopf-Galois extension $B\to P$ are in one-to-one correspondence with
{\em strong connection one forms}.  These are defined as  homomorphisms
$\omega:H\rightarrow\Omega^{1}P$ vanishing on $1$ and satisfying the
following three conditions (see  \cite[Theorem~2.3]{DabGro:str} for
various equivalent descriptions of strong connections):
\begin{subequations}
      \label{eq9}
      \begin{gather}
	\Delta_{\Omega^{1}P}\circ\omega=(\omega\otimes\mId)\circ\mAd ,
	\label{eq9a}\\
	(m_P\otimes\mId)\circ(\mId\otimes\Delta_{R})\circ\omega=
	1\otimes(\mId-\eps) ,
	\label{eq9b}\\
	\mD p-p_{(0)}\omega(p_{(1)})\in(\Omega^{1}B)P,\quad\forall p\in P .
	\label{eq9c}
      \end{gather}
\end{subequations}
A few of the symbols above require explanation. The map
$\Delta_{\Omega^1 P}$
on the left hand side of equation (\ref{eq9a}) is the right caction of
$H$ on $\Omega^1 P$ obtained be restricting the right coaction of $H$ on
the tensor product $P\otimes P$ (right diagonal coaction). Explicitly,
$$
\Delta_{\Omega^1 P}: \Omega^1 P\to \Omega^1 P\otimes H, \quad
\sum_ip^i\otimes \tilde{p}^i\mapsto  \sum_ip^i\sw 0\otimes
\tilde{p}^i\sw 0 \otimes p^i\sw 1 \tilde{p}^i\sw 1.
$$
The map $\mAd$ on the right hand side of equation (\ref{eq9a}) is the
right adjoint coaction of $H$ on itself, i.e., $\mAd: H\to H\otimes H$,
$h\mapsto h\sw 2\otimes S(h\sw 1) h\sw 3$. Thus equation (\ref{eq9a})
corresponds to the classical ad-covariance property of a connection
one-form. The map
$$\chi = (m_P\otimes\mId)\circ(\mId\otimes\Delta_{R})$$
on the left hand side of equation (\ref{eq9b}) has a geometric meaning
of a {\em horizontal lift}, i.e., an operation which lifts an element of
the Lie algebra of a structure group of a principal bundle to a vector
tangent to a fibre. Thus equation (\ref{eq9b}) has the classical
geometric meaning of the property that a connection form evaluated at a
horizontal lift of an element of a Lie algebra, returns back this
element. Finally, equation (\ref{eq9c}) is the strongness condition,
which distinguishes strong connections within a class of all connections
on a Hopf-Galois extension.

Suppose that the antipode $S\in H$ is  invertible.  In this case $P$ is
a left $H$-comodule with the coaction
$$
\Delta_{L}:P\rightarrow H\otimes P, \qquad
	p\mapsto S^{-1}p_{(1)}\otimes p_{(0)}.
$$
One then proves that if
there is a map
$\ell:H\rightarrow P\otimes P$
which satisfies the following conditions
\begin{subequations}
      \label{eq10}
      \begin{gather}
	\ell(1)=1\otimes 1, \label{eq10a}\\
	\chi(\ell(h))=1\otimes h, \label{eq10b}\\
	\ell(h_{(1)})\otimes h_{(2)}=(\mId\otimes\Delta_{R})\circ \ell(h),
	\label{eq10c}\\
	h_{(1)}\otimes \ell(h_{(2)})=(\Delta_{L}\otimes\mId)\circ \ell(h),
	\label{eq10d}
      \end{gather}
\end{subequations}
for all $h\in H$, then the map
\begin{equation}
      \label{eq11}
	\omega:H \rightarrow\Omega^{1}P, \qquad
	h\mapsto \ell(h)-\eps(h)1\otimes 1
\end{equation}
is a strong connection form. This is, in fact, an equivalent description
of strong connection one forms (cf.\ \cite{BrzHaj:rel} for the proof and
discussion in more general situation). Note that equations
(\ref{eq10c})-(\ref{eq10d}) simply state that $\ell$ is an
$H$-bicomodule map, where $H$ is a viewed as a bicomodule via the
coproduct, and $P\otimes P$ is an $H$-bicomodule with coactions
$\Delta_L\otimes P$ and $P\otimes\Delta_R$.

The antipode of the Hopf algebra $H$ of functions on $U(1)$ is
involutive, i.e. $S\circ S =\id$, hence, in particular, invertible. Thus
we can define the left coaction $\Delta_L$ for $P = A^\infty(S^3_\mu)$.
This
again is defined in terms of the $\Z$-grading and comes out as
$$
\Delta_L(x) = u^{-\mdeg (x)}\otimes x,
$$
for any homogeneous $x\in P$. Thus we can follow the above procedure in
the case of   a Hopf-Galois extension $P=A^\infty(S^3_\mu)$  of
$B=A^\infty(S^2_\mu)$,
and define $\ell:H\rightarrow P\otimes P$ by
\begin{subequations}
      \label{eq12}
\begin{align}
          \ell(u^{n})&=\sum_{k=0}^{n}{n\choose k}(a^{\ast})^{n-
k}(b^{\ast})^{k}
	\otimes b^{k}a^{n-k},\label{eq12a}\\
	\ell(u^{-n})&=({1+n\mu})^{-1}\sum_{k=0}^{n}{n\choose k}
	a^{n-k}b^{k}\otimes (b^{\ast})^{k}(a^{\ast})^{n-k},\label{eq12b}\\
	\ell(1)&=1 \otimes 1 , \label{eq12c}
\end{align}
\end{subequations}
for $n\in \N$.
Note that this is simply the expression for
$\can^{-1}(1\otimes h)$
lifted to $P\otimes P$ by omitting the decoration $_B$ on $\otimes_B$
(cf.\ equation (\ref{eq4})).
By definition of $\ell$ the equation (\ref{eq10a}) is satisfied.
Next, since $\chi$ is a lifitng of the canonical map $\can $ to
$P\otimes P$, similar arguments to those used to prove that $\can^{-1}$
is the right inverse of $\can$ ensure that equation (\ref{eq10b}) is
satisfied. Moreover, similarly as in the discussion after the definition
of $\can^{-1}$ by counting the degree it follows
that $\ell$ is a right $H$-colinear map, hence equation (\ref{eq10c})
holds.
Finally since $\mdeg((a^{\ast})^{n-k}(b^{\ast})^{k}) = -n$ and
$\mdeg(a^{n-k}b^{k}) = n$ one easily realises that the map $\ell$ is
also left $H$-colinear, so that  equation (\ref{eq10d}) is satisfied.
Thus we have constructed a strong connection in the quantum contact Hopf
fibration with the strong connection form
$\omega(u^n) = \ell(u^n) - 1\otimes 1$, for all $n\in \Z$.

\subsection{Projectors}

To any Hopf-Galois extension $B\to P$ with the structure Hopf algebra
$H$ and any right $H$-comodule $V$  one can associate a left $B$-module
$\Gamma^H(V,P)$ which plays the role of the module of sections  of the
associated quantum vector bundle. Explicilty, $\Gamma^H(V,P)$ is a
vector space of right $H$-comodule maps $\phi: V\to P$ with the
$B$-action given by $(b\cdot \phi)(v) = b\phi(v)$. If a Hopf-Galois
extension $B\to P$ admits a strong connection and $H$ has a bijective
antipode then for any finite-dimensional $V$, $\Gamma^H(V,P)$ is a
finitely generated projective left $B$-module (cf.\ \cite{DabGro:str}).
Furthermore, the covariant derivative corresponding to a strong
connection gives rise to a connection in module $\Gamma^H(V,P)$.

For any map $\ell:H\rightarrow P\otimes P$ satisfying conditions
(\ref{eq10})  write for all $h\in H$,
\begin{equation}
\label{lsweedler}
\ell(h)=\sum_{i} \ell^{[1]}_{i}(h)\otimes\ell^{[2]}_{i}(h).
\end{equation}
Although the number $r$ of terms in the sum on the right hand side
may depend on $h$, it is always finite. Application of
$\mId\otimes\eps$ to both sides of equation (\ref{eq10b})
yields
\begin{equation}
      \label{eq13}
      \sum_{i}\ell^{[1]}_{i}(h)\ell^{[2]}_{i}(h)=\eps(h).
\end{equation}
Since the left hand side of equation (\ref{eq13}) is finite, for any
$h\in H$ one can define a square matrix $p(h)$ by
\begin{equation}
      \label{eq14}
       P\ni p(h)_{ij}=\ell^{[2]}_{i}(h)\ell^{[1]}_{j}(h).
\end{equation}
Equation (\ref{eq13}) immediately implies  that
$\sum_{k}p(h)_{ik}p(h)_{kj}=\eps(h)p(h)_{ij},$
hence the matrix $p(h)$ is an idempotent, provided $\eps(h)=1$.
Furthermore, if $h$ is a grouplike element, i.e., $\Delta(h) = h\otimes
h$, then equations (\ref{eq10c}), (\ref{eq10d}) imply that
$\Delta_R(p(h)_{ij}) = p(h)_{ij}\otimes 1$, i.e., $p(h)_{ij} \in B$.
In other words, for any grouplike element $h\in H$,
$p(h)$ is an $r\times r$ matrix
($r$ being the number of terms in (\ref{lsweedler})),
which is an idempotent in a matrix ring over $B$. Therefore it defines
a finitely generated projective module via $B^r\, p(h)$.
The connection corresponding to $\ell$ is simply the Grassmann
or Levi-Civita connection in $B^{r}\, p(h)$,
i.e., a connection determined by the idempotent (cf.\ \cite{CunQui:alg}).

In the case of the Hopf-Galois extension $A^\infty(S^2_\mu)\to
A^\infty(S^3_\mu)$
each of the $u^n$ is a grouplike element. Therefore  the map $\ell$
defined by formulae (\ref{eq12})  gives an infinite family of
idempotents $p(u^{n})$, $n\in\mZ$.  Each of the $p(u^n)$ is an
$(n+1)\times (n+1)$-matrix with entries from $A^\infty(S^2_\mu)$ (the
latter
claim can be easily confirmed by the degree counting). Obviously there
is an ambiguity in factorising $\ell$ into
$\ell^{[1]}\otimes\ell^{[2]}$ (scalar coefficients can be factorised
in infinitely many ways into legs of tensor product). However if one
requires
$p(u^{n})$ to be Hermitian (i.e. projectors in $B$) then the
unique possibility turns out to be
\begin{subequations}
      \label{eq15}
\begin{align}
	\ell(1)&=1\otimes 1,
\nonumber\\
          \ell(u^{n})&=\sum_{k=0}^{n}
	\left[\sqrt{{n\choose k}}(a^{\ast})^{n-k}(b^{\ast})^{k}\right]
	\otimes
	\left[\sqrt{{n\choose k}}b^{k}a^{n-k}\right],
\\
	\ell(u^{-n})&=
	\sum_{k=0}^{n}\left[\sqrt{{n\choose k}}({1+n\mu})^{-1}
        a^{n-k}b^{k}\right] \otimes
	\left[\sqrt{{n\choose k}}(b^{\ast})^{k}(a^{\ast})^{n-k}\right],
\end{align}
\end{subequations}
for $n\in \N$.
This choice leads to an infinite family of Hermitian
projectors with entries from
$A^\infty(S_{\mu}^{2})$. Explicitly, $p(1)= 1$, and
\begin{subequations}
      \label{proj}
\begin{align}
         p(u^{n})_{kl}&=
	\sqrt{{n\choose k}{n\choose l}}b^{k}a^{n-k}(a^{\ast})^{n-
l}(b^{\ast})^{l},
	 \label{proja} \\
	p(u^{-n})_{kl}&=
	\sqrt{{n\choose k}{n\choose l}}(b^{\ast})^{k}(a^{\ast})^{n-
k}({1+n\mu})^{-1}a^{n-l}b^{l},
	\label{projb}
\end{align}
\end{subequations}
$n\in \N$.
At this point it is interesting to mention that apparently
these fromulae are polynomial in $a$, $b$, $a^*$, $b^*$,
and polynomial only in $({1+n\mu})^{-1}$, with $n\in \N$,
but not in $\mu$ (e.g. (\ref{projb})).
However, their entries have to be properly rearranged
in order to express them
in terms of the generators $X$, $Z$, $Z^{\ast}$.
In $A^\infty(S^{3}_{\mu})$ this
can be always done with the help of relations (\ref{eq1}) and (\ref{invs})
at the cost of creating new expressions $({1+n\mu})^{-1}$
in (\ref{proja}) and $({1-n\mu})^{-1}$ in (\ref{projb}),
with $n\in \N$. Alltogether the whole set of projectors can be rewritten
in terms of $X$, $Z$, $Z^{\ast}$ and all $({1+k\mu})^{-1}$ with $k\in \Z$.
(The reason for this behaviour is that the the elements
$X$, $Z$, $Z^{\ast}$ and $\mu$ do not generate the whole grade-zero
polynomial subalgebra of $A(S^3_\mu) $.)
For instance the first few projectors
$p(u)$, $p(u^{-1})$, $p(u^{2})$, $p(u^{-2})$
come out in a matrix form as
\begin{subequations}
   \label{projmatrix}
      \begin{align}
	p(u)&=\left(
	\begin{array}{cc}
	    \frac{1}{2}(1+ \mu)+X & Z\\
	    Z^{\ast} & \frac{1}{2}(1+ \mu)-X
	\end{array}
	\right),
\\
~\nonumber\\
	p(u^{-1})&=\left(
	\begin{array}{cc}
	     \frac{1}{2}(1- \mu)+X & Z^{\ast}\\
	    Z &  \frac{1}{2}(1- \mu)-X
	\end{array}
	\right),
\\
~\nonumber\\
	p(u^{2})&=\frac{1}{1+\mu}
\\
         &\times\left(
	\begin{array}{ccc}
	    \left(X+\frac{1+\mu}{2}\right)\left(X+\frac{1+3\mu}{2}\right)&
	    \sqrt{2}\left(X+\frac{1+3\mu}{2}\right)Z&
	    Z^{2}\\
	    &&\\
	    \sqrt{2}Z^{\ast}\left(X+\frac{1+3\mu}{2}\right)&
	    2 \left(\frac{1+\mu}{2}+X\right)\left(\frac{1+\mu}{2}-X\right) &	    \sqrt{2}\left(\frac{1+\mu}{2}-X\right)Z\\
	    &&\\
	    (Z^{\ast})^{2}&
	    \sqrt{2}Z^{\ast}\left(\frac{1+\mu}{2}-X\right)&
	    \left(\frac{1+\mu}{2}-X\right)\left(\frac{1+3\mu}{2}-X\right)
	\end{array}
	\right),
\nonumber
\\
~\nonumber\\
	p(u^{-2})&=\frac{1}{1-\mu}
\\
         &\times\left(
	\begin{array}{ccc}
	    \left(X+\frac{1-\mu}{2}\right)\left(X+\frac{1-3\mu}{2}\right)&
	    \sqrt{2}\left(X+\frac{1-3\mu}{2}\right)Z^{\ast}&
	    (Z^{\ast})^{2}\\
	    &&\\
	    \sqrt{2}Z\left(X+\frac{1-3\mu}{2}\right)&
	    2 \left(\frac{1-\mu}{2}+X\right)\left(\frac{1-\mu}{2}-X\right) &	    \sqrt{2}\left(\frac{1-\mu}{2}-X\right)Z^{\ast}\\
	    &&\\
	    Z^{2}&
	    \sqrt{2}Z\left(\frac{1-\mu}{2}-X\right)&
	    \left(\frac{1-\mu}{2}-X\right)\left(\frac{1-3\mu}{2}-X\right)
	\end{array}
	\right).
 \nonumber
      \end{align}
\end{subequations}

Note an interesting symmetry
between $p(u^{n})$, and $p(u^{-n})$ for low values of $n$.
The projector  $p(u^{-n})$  is obtained from the projector $p(u^{n})$ by
replacing $\mu$ by $-\mu$ and interchanging of $Z$ with $Z^*$.
This is true for any value of charge $n$ as can be verified directly
from the explicit expressions for greater charges $n$,
which can be presented.
However this follows also from the following symmetry properties.
First observe that the transformation
\begin{equation}
Z\mapsto Z^*, \quad Z^*\mapsto Z, \quad X\mapsto X, \quad \mu\mapsto -\mu
\label{eq.auto}
\end{equation}
does not affect the defining relations (\ref{eq3})
and defines an automorphism of the algebra $A^\infty(S^{2}_{\mu})$.
This symmetry of $A^\infty(S^{2}_{\mu})$ comes in fact from the
following symmetry
of $A^\infty(S^{3}_{\mu})$.
Using the elements $\sqrt{1+k\mu}\in A^\infty(S^{3}_{\mu})$,
their inverses $1/\sqrt{1+k\mu}$ for $k\in \Z$ and the relations
(\ref{sqrts}) we see that the map
\begin{equation}
  \label{vartheta}
\vartheta: A^{\infty}(S^{3}_{\mu})\to A^\infty(S^{3}_{\mu}), \qquad
\mu\mapsto -\mu, \quad a\mapsto A=\sqrt{1-\mu}a^{\ast}, \quad
b\mapsto B=\sqrt{1-\mu}b^{\ast},
\end{equation}
extends to an algebra automorphism. Note that $\vartheta$ maps degree $n$
elements to degree $-n$ elements, and on the level of degree $0$
elements corresponds to the automorphism of  $A^\infty(S^{2}_{\mu})$ given by
(\ref{eq.auto}).
Now, using the relations
$$
      a^{n}(1-\mu)=\frac{1+(n-1)\mu}{1+n\mu}a^{n}, \quad
      b^{n}(1-\mu)=\frac{1+(n-1)\mu}{1+n\mu}b^{n}, \quad n>0
$$
and the fact that $\mu$ is central in $A^\infty(S^{2}_{\mu})$ we find
that
$$
\vartheta(p(u^n)_{kl}) = p(u^{-n})_{kl}.
$$
Thus the authomorphism $\vartheta |_{A^\infty(S^{2}_{\mu})} : A^\infty(S^{2}_{\mu})\to
A^\infty(S^{2}_{\mu})$  turns the degree $n$ projectors into degree $-n$
projectors, and its existence proves the stated symmetry of monopole
projectors.

\section{Final remarks}

If $\mu$ is not required to be invertibile
the underlying polynomial $*$-algebra $A(S^3_\mu)$
admits a specification to $\mu = 0$, after which it coincides
with the usual $*$-algebra of polynomials on the `classical' $S^3$.
As far as the polynomial $*$-algebra $A'(S^2_\mu)$ is concerned
since $\mu$ is central it can be specified to any nonzero real value.
This yields a family of quantum 2-spheres isomorphic to
the universal enveloping algebra of $su(2)$
with a constrained value $1/\mu^2$
of the quadratic Casimir element (c.f. \cite{Dab:gar}).
The $*$-algebra $A(S^2_\mu)$ (i.e., when the invertibility of $\mu$ is
not assumed) admits in addition a specification to $\mu = 0$,
which clearly corresponds to polynomials on the `classical' $S^2$.
Note that the equations (\ref{proj}) with $\mu$ considered
as a real deformation parameter rather than as a central generator
define also a family of projectors $p(u^n)$ over such quantum 2-spheres.
In particular, when $\mu=0$, the projectors $p(u^n)$ correspond to line
bundle projectors of the monopole charge or Chern's number $n$ over $S^2$.
Note also that then $\vartheta$ defined by (\ref{vartheta})
is the orientation reversing automorphism of $S^3$.

It would be interesting to investigate further our bundles
and the associated projective modules, especially the Chern-Connes
pairing between $K$-theory and $K$-homology.
For that aim it would be convenient to have a $C^*$-algebraic
version of the construction above as otherwise computation of e.g.\ the
$K$-groups is a formidable task. It can be seen that the $*$-algebras
$A(S^3_\mu)$ and $A(S^2_\mu)$ admit certain $C^*$-algebraic completions.
Although intuitively resembling some topological quantum four
(resp. three) dimensional spaces rather than 3-spheres (resp. 2-spheres),
they nevertheless would constitute interesting examples.
Unfortunately, there is one obstacle for this task.
In order to write our formulae (\ref{eq4})
for the inverse of the canonical map and then for the connection
and the projectors we should adjoin to
$A(S^3_\mu)$ and to $A(S^2_\mu)$  an infinite number of elements
$({1+k\mu})^{-1}$, $k\in \Z$.
(To be able to implement the symmetry discussed at the end of Section 3.2 we
should adjoin additionally the elements $\sqrt{1+k\mu}$ and $1/\sqrt{1+k\mu}$.)
It can be seen that this spoils not only the $C^*$-algebraic completion
but even the $*$-algebraic version. This follows from
the representation theory of e.g.\ $A(S^2_\mu)$,
which can be inferred from that of $su(2)$.
In fact all bounded representations decompose into finite dimensional
irreducible ones. These in turn are as follows.
There is a family of one-dimensional representations (characters)
parametrised by the points of $S^2$, which represent $\mu$ by $0 $
(they obviously do not extend to $A'(S^2_\mu)$).
In addition, in each dimension $N\in \N$ there are two $*$-representations,
labeled by $\sigma = \pm 1$, which represent $\mu$ either
by $1/N $ or by $-1/N$.
Namely, they are given by
\begin{equation}
\begin{array}{l}
 \mu ~v_{m}  = \frac{\sigma}{N} ~v_{m},\\
 X ~v_{m}  = \frac{\sigma m}{N} ~v_{m},\\
 Z ~v_{m} = \frac{\sigma}{2N}
\sqrt{({N+1}-2m)({N-1}+2m)} ~v_{m-1},\\
 Z^* ~v_{m} = \frac{\sigma}{2N}
\sqrt{({N-1}-2m)({N+1}+2m)} ~v_{m+1},
\end{array}
\end{equation}
with respect to an orthonormal basis $v_{m}$, where
$m\in\{-\frac{N-1}{2}, -\frac{N-3}{2}, \dots , \frac{N-3}{2}, \frac{N-1}{2}\}$.

This excludes a possibility of adjoining the needed elements
if we want our $C^*$-algebra to describe more than merely the commutative $S^2$.
At this point it is interesting to make the following observation.
Had we constrained ourselves to construct just certain subclass of projectors,
say for a selected class of charges,
we might have a chance to accomplish a nondegenerate $*$-algebra
extended by some (but not all) elements $({1+k\mu})^{-1}$
and then also its $C^*$-algebraic version.
For instance restricting only to the positive charges $n$ leaves
at our disposal the series of all representations that
represent $\mu$ by $1/N $.
It is not clear however what could be a possible interpretation
of such a breaking of the (magnetic) charge conjugation
by the contact  structure quantisation with a noncentral parmeter.

\section*{Acknowledgements}
Tomasz Brzezi\'nski  would like to thank the Engineering and Physical
Sciences Research
Council for an Advanced Fellowship. The research of  Bartosz Zieli\'nski
is supported by the EPSRC grant GR/S01078/01

\end{document}